\documentclass[12pt]{article}%%25.11.2016
\usepackage{amsmath}
\usepackage{amssymb}
\usepackage{dsfont}
\usepackage{cite}
\newsymbol \blackbox 1004
\newcommand{\eh}{\hfill}\newlength{\sperr}

\newenvironment{proof}{{\settowidth{\sperr}{\bf\rm
Proof}%
\par\addvspace{0.3cm}\noindent\parbox[t]{1.3\sperr}
{\bf\rm P\eh r\eh o\eh o\eh f\eh }%
}}{\nopagebreak\mbox{}
$\blackbox$\par\addvspace{0.3cm}}

\def\cls{\mathcal{S}}

\def\nn{\nonumber}

\def\Lam{\Lambda}
\def\s{\sigma}
\def\la{\lambda}

\def\wt{\widetilde}
\def\ov{\overline}

\def\p{\partial}

\def\BC{{\mathbb C}}
\def\BR{{\mathbb R}}
\def\BN{{\mathbb N}}

\def\cla{{\mathcal A}}

\def\cli{{\mathcal I}}
\def\clj{{\mathcal J}}
\def\cll{{\mathcal L}}

\def\diag{\mathrm{diag}}
\newcommand{\E}{\mathrm{e}}
\newcommand{\I}{\mathrm{i}}

\newtheorem{Pa}{Paper}[section]
\newtheorem{Tm}[Pa]{{\bf Theorem}}

\newtheorem{Cy}[Pa]{{\bf Corollary}}
\newtheorem{Rk}[Pa]{{\bf Remark}}

\newtheorem{Pn}[Pa]{{\bf Proposition}}

\title{Continuous and discrete dynamical Schr\"odinger systems: explicit solutions}

\author{B. Fritzsche, B. Kirstein, I.Ya. Roitberg, A.L. Sakhnovich}

\date{}
\begin{document}
\maketitle

\begin{abstract} We consider continuous and discrete Schr\"odinger systems with self-adjoint matrix potentials
and with additional dependence on time (i.e.,  dynamical Schr\"odinger systems). 
Transformed and explicit solutions are constructed using our generalized (GBDT) version of the B\"acklund-Darboux transformation.
Asymptotic expansions of these solutions in time are of interest.

\end{abstract}

{MSC(2010):   35Q41, 37C80, 39A12, 47B36.} 

Keywords:  {\it  Schr\"odinger equation, dynamical system, Jacobi matrix, \\  B\"acklund-Darboux transformation, dispersion, explicit solution.}
%rectangular matrix potential, pseudo-exponential potential, direct problem, inverse problem, 
%rational matrix function, 
%realization.} 

\section{Introduction} \label{intro}
\setcounter{equation}{0} 
Dynamical Dirac and Schr\"odinger systems play an essential role in mathematical physics and are actively studied, especially in the recent years
(see, e.g., \cite{AMRa, BeMi, BoPP, EKMT, EKT, KT, PrdO, SaA15, SaA17, T1} and numerous references therein).
Continuous dynamical Schr\"odinger system has the form:
\begin{align}\label{I2}&
\I \frac{\p}{\p t}\psi(x,t)=\big(H \psi\big)(x,t), \quad H:=-\frac{\p^2}{\p x^2}+u(x) \quad (u=u^*),
\end{align}
where $u$ is an $h \times h$ matrix function, $h\in \BN$, and $\BN$ is the set of natural numbers.
The  matrix function $u$ is called the {\it potential} of \eqref{I2} and this potential does not depend on $t$ in our case.

In discrete dynamical Schr\"odinger system we use Jacobi matrices $\clj$ instead of $H$
since Jacobi operators  ``can be viewed as
the discrete analogue of Sturm-Liouville operators'' \cite[Preface]{T0}. The corresponding system is given by the
formula:
\begin{align}\label{I3}&
\I \left(\frac{\p}{\p t}\Psi\right)(t)=\clj \Psi(t), 
\end{align}
where $\clj$ is a semi-infinite block Jacobi matrix and $\Psi$ is a block vector 
\begin{align}\label{I4}&
\clj=\begin{bmatrix} b_1 & a_1 & 0 & 0 & 0 & \ldots \\ 
c_2 &b_2 & a_2  & 0 & 0 & \ldots \\
0 & c_3 &b_3 & a_3   & 0 & \ldots \\
 \ldots & \ldots & \ldots & \ldots & \ldots & \ldots 
\end{bmatrix}, \quad \Psi(t)=\begin{bmatrix} \psi_1(t) \\ \psi_2(t) \\ \psi_3(t)\\ \ldots \end{bmatrix}.
\end{align}
Here, the blocks $a_k$, $b_k$ and $c_k$ are $h \times h$ matrices and $c_{k+1}=a_{k}^*$ $(k\geq 1)$. 

Explicit solutions of dynamical systems are important as models and examples and they are
also essential in applications. Various explicit solutions of time-independent systems were 
constructed using commutation methods \cite{D, Ge, GeT, KoSaTe} and several versions of B\"acklund-Darboux transformations.
B\"acklund-Darboux transformation is a well-known tool in the spectral theory and theory of  explicit solutions.
The original equation, which was studied by Darboux, is the Schr\"odinger equation
\begin{align}\label{h1}
-y^{\prime \prime}(x, \la)+u(x)y(x,\la)=\la y(x,\la)\quad (u=u^*), \quad y^{\prime}:=\frac{d}{d x}y.
\end{align}
Later, and especially in the last
40 years, this transformation was greatly modified, generalized and applied to a variety of linear and nonlinear
equations (see, e.g., \cite{Ci, Gu, Mar, MS, SaSaR}). 

It was shown recently  in \cite{SaA17, SaA16} that the GBDT version of B\"acklund-Darboux transformation
(for GBDT see \cite{SaA2, SaA3, ALS-PT, ALS-KP, SaA07, SaSaR} and references therein) may be successfully applied to the construction of explicit solutions of dynamical systems
as well. 

In the present paper, we consider the important case of  continuous dynamical Schr\"odinger system and a more difficult case of discrete system
(i.e., system \eqref{I3}). Some preliminaries are presented in Section \ref{GBDT}, continuous dynamical Schr\"odinger system is dealt with in Section \ref{SchrEq}
and system \eqref{I3} is considered in Section \ref{DscrSchrEq}.

 The dependence of our solutions of \eqref{I2} and \eqref{I3} on time is described by the factor $\E^{\I t A}$,
where $A$ is a parameter matrix (generalized eigenvalue) of the GBDT transformation. Since $A$ is not necessarily self-adjoint and may have Jordan cells of different orders,
the asymptotic expansion of our solutions of \eqref{I2} and \eqref{I3}  essentially differs (see  Remark \ref{DD}) from the classical Jensen-Kato formulas (see
\cite{JK} as well as  further references in \cite{EKMT, EKT}).

As usual, $\BR$ denotes the set of real values,  $\BN$ is the set of natural numbers, and the complex plane is denoted by $\BC$.
Notation $\cls^*$ stands for the matrix which is the conjugate transpose of $\cls$, we write $\cls >0$ when $\cls$ is a positive-definite matrix, and
$I_m$ stands for  the $m\times m$ identity matrix. The notation $J=\diag\{J_1, J_2, \ldots \}$ means that $J$ is a diagonal or block diagonal matrix
with the entries (or block entries) $J_1$, $J_2$ and so on.

%%%%%%%%%%%%%%%%%%%%%%%%%%%%%%%%%%%%%%%%%%%%%%
%%%%%%%%%%%%%%%%%%%%%%%%%%%%%%%%%%%%%%%%%%%%%%%
\section{GBDT: preliminaries} \label{GBDT}
\setcounter{equation}{0}
GBDT (generalized B\"acklund-Darboux transformation) was first introduced in \cite{SaA2}, and 
a more general version of GBDT for first order systems rationally depending on the spectral parameter
(in particular, for systems of the form $w^{\prime}= G(x, \lambda ) w$, $G(x, \lambda )=- \sum_{k=-r}^{r}
\lambda^{k}q_{k}(x)$)
was treated in \cite{SaA3, SaSaR} (see also some references therein). First order system 
\begin{equation} \label{bt1}
w^{\prime}(x, \lambda)= G(x, \lambda ) w(x, \lambda ),
\hspace{1em}
G(x,\la)=-\la q_1(x)-q_0(x),
\end{equation}
where $w$ takes values in $\BC^m$ $(m:=2h)$ and $m\times m$ coefficients $q_1$ and $q_2$ have the form
\begin{equation} \label{bt11}
q_{1}= \left[  \begin{array}{lr} 0 & 0 \\ I_{h} &  0 \end{array}
\right], \hspace{1em} q_{0}(x)=- \left[  \begin{array}{lr} 0 &
I_{h} \\  u(x) &  0 \end{array} \right], \quad u(x)=u(x)^*,
\end{equation}
is equivalent to the matrix Schr\"odinger equation \eqref{h1} with a self-adjoint $h \times h$ potential $u(x)$.
Here we present basic GBDT  results for  this system (see, e.g., \cite{SaA3, ALS-PT}). The connection
with the Schr\"odinger equations \eqref{I2} and \eqref{h1}  is discussed in greater detail in the next section.
\begin{Rk} \label{RkD1} We consider systems \eqref{bt1} and \eqref{I2} on finite or infinite intervals $\cli$, that is, we assume
that $x\in \cli$. Without loss of generality we assume also that $0 \in \cli$ and speak later about parameter
matrices $S(0)$ and $\Pi(0)$  instead of $S(x_0)$ and $\Pi(x_0)$ for some fixed $x_0 \in \cli$.
The most interesting for us is the case of the semiaxis $\cli=[0,\infty)$.
\end{Rk}

In general, GBDT is determined by the choice of  5
parameter matrices (this case was treated in \cite{ALS-PT}, where $u$ was not
necessarily self-adjoint). However, relations \eqref{bt11} (including $u=u^*$) imply additional equalities:
\begin{equation} \label{h4}
q_k^*=-jq_k j^* \quad (k=0,1), \quad j:=\begin{bmatrix} 0 & I_h \\ - I_h &0
\end{bmatrix}, \quad j^*=j^{-1}=-j.
\end{equation}
Thus the conditions of  Proposition 1.4 from \cite{SaA3} are fulfilled, and we may use this proposition
and some formulas from its proof. The following
statements in this section are particular cases of \cite[Theorem 2.1]{ALS-PT} (or \cite[Theorem 1.2]{SaA3})
completed by \cite[Proposition 1.4]{SaA3}.

Hence, in the present case we use 3 parameter matrices.
More precisely, we choose some {\it initial} system \eqref{bt1} (or, equivalently, the initial potential $u=u^*$ of Schr\"odinger equation \eqref{h1}) and fix
$n \in \BN$. Then, we fix 
$n\times n$ matrices $A$ and $S(0)=S(0)^*$, and an $n\times m$ 
$(m=2h)$ matrix $\Pi(0)$ such that
the following matrix identity holds:
\begin{equation} \label{bt2}
AS(0)-S(0)A^*= \Pi(0)j \Pi(0)^{*}.
\end{equation}
Suppose that such parameter matrices are
fixed and that the potential $u(x)$ 
is locally summable on $\BR$.
Now, we can introduce matrix functions $\Pi(x)$
 and $S(x)$ with the values $ \Pi(0)$ and $S(0)$ at $x=0$
as the solutions of the linear
differential equations
\begin{equation} \label{h2}
\Pi^{\prime}= A \Pi q_{1}+\Pi q_{0},
\quad  S^{\prime}=  \Pi q_{1}j
\Pi ^{*} ,
\end{equation}
where $q_1$ and $q_0$ are given by \eqref{bt11}, and so $q_{1}j=(q_{1}j)^*$. Thus, in view of $S(0)=S(0)^*$,
we have
\begin{equation} \label{h10}
S(x)=S(x)^*.
\end{equation}
Notice that equations (\ref{h2}) are constructed in such a way
that the identity
\begin{equation} \label{bt4}
A S(x)-S(x)A^*= \Pi (x)j \Pi (x)^{*}
\end{equation}
follows (for all $x\in \BR$) from (\ref{bt2}) and (\ref{h2}).  (The relation is obtained by
the direct differentiation
of  the both sides of (\ref{bt4}).) Assuming that $\det
S(x) \not \equiv 0$ we can define a matrix function
\begin{equation} \label{bt5}
w_{A}(x, \lambda )=I_{m}- j\Pi(x)^{*}S(x)^{-1}(A- \lambda
I_{n})^{-1} \Pi(x),
\end{equation}
where $\la \not\in \s(A)$ ($\s$ means spectrum).
\begin{Tm} \label{Tmbt} 
Suppose that   the relation  \eqref{bt2}
is valid, and matrix functions  $\Pi (x)$  and $S(x)$
satisfy  equations \eqref{h2}
where \eqref{bt11} holds. Then, in the points of
invertibility of $S(x)$, the matrix function $w_A(x,\la)$ satisfies the
system
\begin{align} & \label{bt8}
w_A^{\prime}(x, \lambda )= \wt{G}(x, \lambda ) w_A(x, \lambda )
-w_A(x, \la)G(x, \la), \\
& \label{h3}
 \wt{G}(x, \lambda ):= -\la  q_1(x)-\wt q_0(x),
\end{align}
where  the
coefficient $\wt{q}_{0}(x)$ is given by the formula
\begin{align} & \label{bt9}
\wt{q}_{0} =q_{0} - (q_{1} j X
-j X q_{1}), \quad X(x): = \Pi(x)^{*}S(x)^{-1}\Pi(x).
\end{align}
\end{Tm}
\begin{Rk} \label{RkSpos1} Formulas \eqref{bt11} and  \eqref{h2} yield 
\begin{align} & \label{bt9!}
q_1 j=\begin{bmatrix}0 & 0 \\ 0 & I_h\end{bmatrix}\geq 0, \quad S^{\prime}(x)\geq 0,
\end{align}
and so  $S(x)>0$  for $x\geq 0$ under additional condition  $S(0)>0$.
In particular, the condition of invertibility of $S(x)$ from  Theorem \ref{Tmbt}
is fulfilled automatically when $\cli=[0,\, \infty)$ and $S(0)>0$.

The matrix functions $S(x)^{-1}$, $X(x)$ and $w_A(x,\la)$ are well-defined
in  this case.
\end{Rk}
According to Theorem \ref{Tmbt}, the multiplication by $w_A$
transforms each solution $w$ of (\ref{bt1}) into the
 solution $\widetilde{w}=w_Aw$ of the system
$\widetilde{w}^{\prime}= \widetilde{G} \widetilde{w}$ with the
coefficients  of $\wt{G}$ given by \eqref{h3} and (\ref{bt9}). This
transformation of the  solutions $w$ and coefficients
$q_k$ is called GBDT. Matrix function $w_A$ is the so called
Darboux matrix. The right hand side of \eqref{bt5}
(with  the additional property \eqref{bt4} and $x$ fixed)
has the form of the Lev Sakhnovich transfer matrix function
\cite{SaSaR, SaL1, SaL2}.

Under the conditions of Theorem \ref{Tmbt} we have also
\begin{equation} \label{h8}
( \Pi^{*}S^{-1})^{\prime}(x)= {q}_{1}^*\Pi(x)^{*}S(x)^{-1}A+\wt q_0(x)^*
\Pi(x)^{*}S(x)^{-1}.
\end{equation}
Clearly, the definition \eqref{bt9} of $X$ and formula \eqref{h10} imply that
\begin{align} & \label{h11}
X(x)=X(x)^*.
\end{align}
%\begin{equation} \label{bt10.2'}
%(S^{-1} \Pi_{1})^{\prime}(x)= \sum_{p=0}^{r}A_{2}^{p}S(x)^{-1}
%\Pi_{1}(x)
%\wt{q}_{p}(x).
%\end{equation}

%%%%%%%%%%%%%%%%%%%%%%%%%%%%%%%%%%%%%%%%%%%%%%%
%%%%%%%%%%%%%%%%%%%%%%%%%%%%%%%%%%%%%%%%%%%%%
%%%%%%%%%%%%%%%%%%%%%%%%%%%%%%%%%%%%%%%%%%%%%%
\section{Explicit solutions of the dynamical system \eqref{I2} and GBDT of the matrix Schr\"odinger equation} \label{SchrEq}
\setcounter{equation}{0}
\paragraph{1.}
Let us write down the coefficient $\wt q_0$ of the transformed system in the block form.
We partition $\Pi$ into two $h \times h$ blocks and partition $X$ introduced in \eqref{bt9} into
four $h \times h$ blocks:
\begin{align} & \label{h5}
\Pi= \begin{bmatrix} \Phi_1 & \Phi_2 \end{bmatrix}, \quad X=\{X_{ij}\}_{i,j=1}^2, \quad X_{ij}=\Phi_i^* S^{-1}\Phi_j .
\end{align}
Thus, $\wt q_0$ in Theorem \ref{Tmbt} (see \eqref{bt9}) has the form
\begin{align} & \label{h6}
\wt q_0= -\begin{bmatrix} - X_{22} & I_h \\ u+X_{12}+X_{21} & X_{22} \end{bmatrix}.
\end{align}
In order to rewrite \eqref{h8} in a more convenient form, we shall need also the block representation of $\Pi^*S^{-1}$,
$q_1^*$ and $\wt q_0^*$:
\begin{align} & \label{D1}
Z=\begin{bmatrix} z_1 \\ z_2 \end{bmatrix}:=\Pi^*S^{-1}, \quad q_1^*=\left[  \begin{array}{lr} 0 &  I_{h}\\ 0 &  0 \end{array}
\right],
 \\ & \label{D2}
 \wt q_0^*= -\begin{bmatrix} - X_{22} &  u+X_{12}+X_{21}\\ I_h  & X_{22} \end{bmatrix},
\end{align}
which follows from \eqref{bt11}, \eqref{h11} and \eqref{h6}. Now, \eqref{h8} takes the form
\begin{align} & \label{D3}
z_1^{\prime}=z_2A+X_{22}z_1-(u+X_{12}+X_{21})z_2, \quad z_2^{\prime}=-z_1-X_{22}z_2.
\end{align}
Differentiating the second equality in \eqref{D3} (and taking into account the first equality), we obtain
\begin{align} & \label{D4}
z_2^{\prime \prime}=-z_2A+(u+X_{12}+X_{21}-X_{22}^{\prime}+X_{22}^2)z_2.
\end{align}
Using \eqref{D4} we derive the main theorem in this section
\begin{Tm} \label{TmCS}
Let the parameter matrices $A$, $S(0)=S(0)^*$ and $\Pi(0)$ be chosen so that  \eqref{bt2}
is valid, let the $h \times h$ potential $u=u^*$ be  locally summable on $\BR$, and introduce  $\Pi (x)$  and $S(x)$
via \eqref{h2}
where \eqref{bt11} holds. 

Then, in the points of
invertibility of $S(x)$, the matrix function 
\begin{align} & \label{D5}
\psi(x,t)=\begin{bmatrix} 0 & I_h \end{bmatrix}\Pi(x)^*S(x)^{-1}\E^{-\I t A}
\end{align}
satisfies the continuous dynamical Schr\"odinger system
\begin{align}\label{I2'}&
\I \frac{\p}{\p t}\psi(x,t)=\big(\wt H \psi\big)(x,t), \quad \wt H:=-\frac{\p^2}{\p x^2}+\wt u(x),
\end{align}
where $\wt u =\wt u^*$ is given by the formula
\begin{equation} \label{bt15}
\widetilde{u}(x)=u(x)+2\big(X_{12}(x)+X_{21}(x)+X_{22}(x)^2\big).
\end{equation}
\end{Tm}
\begin{proof}. Taking into account \eqref{h2}, \eqref{h8} and \eqref{bt11},  \eqref{D2}, we calculate  $X_{22}^{\prime}$:
\begin{align} \nn 
X_{22}^{\prime} &=\begin{bmatrix} 0 & I_h \end{bmatrix}(\Pi^*S^{-1})^{\prime}\Pi\begin{bmatrix} 0 \\  I_h \end{bmatrix}+
\begin{bmatrix} 0 & I_h \end{bmatrix}\Pi^*S^{-1}\Pi^{\prime} \begin{bmatrix} 0 \\  I_h \end{bmatrix}
\\ \nn &
= \begin{bmatrix} 0 & I_h \end{bmatrix}\wt q_0^* \Pi^*S^{-1}\Pi\begin{bmatrix} 0 \\  I_h \end{bmatrix}-
\begin{bmatrix} 0 & I_h \end{bmatrix}\Pi^*S^{-1}\Pi  \begin{bmatrix} I_h \\  0 \end{bmatrix}
\\ & \label{h9}
=- X_{12}-X_{22}^2-X_{21}.
\end{align}
In view of \eqref{h9} (and definition \eqref{bt15} of $\wt u$), we rewrite \eqref{D4} in the form
\begin{align} \nn
z_2(x)A &=-z_2^{\prime \prime}(x)+\big(u(x)+2\big(X_{12}(x)+X_{21}(x)+X_{22}(x)^2\big)\big)z_2(x)
\\  & \label{D6}
=-z_2^{\prime \prime}(x)+\widetilde{u}(x)z_2(x).
\end{align}
According to \eqref{D1}, we have $z_2=\begin{bmatrix} 0 & I_h \end{bmatrix}\Pi^*S^{-1}$. Therefore,
\eqref{D5} and \eqref{D6} imply \eqref{I2'}. We also note that  $\wt u =\wt u^*$ is immediate from
$u=u^*$ and formulas \eqref{h11} and \eqref{bt15}.
\end{proof}
\begin{Rk} \label{RkSpos} Under conditions  $S(0)>0$ and  $x \in [0,\, \infty)$, 
the matrix function $S(x)$ is invertible $($recall Remark \ref{RkSpos1}$)$. Thus, the matrix functions
$S(x)^{-1}$, $X(x)$, $\wt u(x)$ and $\psi(x)$
 $($considered in Theorem  \ref{TmCS}$)$
 are well-defined under  these conditions.
\end{Rk}
%%%%%%%%%%%%%%%%%%%%%%%%%%%%%%%%%%%%%
\paragraph{2.} 
If the conditions of Theorem \ref{TmCS} and Remark \ref{RkSpos} are valid, we obtain
the following corollary.
\begin{Cy} \label{CyDD} Consider dynamical Schr\"odinger equations on $\cli =[0,\infty)$, let the conditions of Theorem \ref{TmCS} hold, 
and assume that $S(0)>0$. Then, the columns of the matrix function $\begin{bmatrix} 0 & I_h \end{bmatrix}\Pi(x)^*S(x)^{-1}$
belong to $L^2_h(0,\infty)$ $($i.e., these columns are squarely summable$)$ and the solutions $\psi(x,t)g$ $(g \in \BC^h)$ of system \eqref{I2'} belong
to $L^2_h(0,\infty)$ for each fixed $t$.
\end{Cy}
\begin{proof}. In view of the second equality in \eqref{h2} and the first equality in \eqref{bt9!}, we have
\begin{align}  & \label{DD1}
\big(S(x)^{-1}\big)^{\prime}=-S(x)^{-1}\Pi(x)\begin{bmatrix}0 & 0 \\ 0 & I_h\end{bmatrix}\Pi(x)^*S(x)^{-1}.
\end{align}
Formula \eqref{DD1} implies that
\begin{align}  & \nn
\int_0^{\ell} \Big(\begin{bmatrix} 0 & I_h \end{bmatrix}\Pi(x)^*S(x)^{-1}\Big)^*\begin{bmatrix} 0 & I_h \end{bmatrix}\Pi(x)^*S(x)^{-1}dx
\\ & \label{DD2} \,
=S(0)^{-1}-S(\ell)^{-1}<S(0)^{-1},
\end{align}
which proves the corollary.
\end{proof}
\begin{Rk} \label{DD}
For the study of the dependence on time of the solutions \eqref{D2} and \eqref{J2} of the discrete and continuous, respectively, Schr\"odinger equations, 
one may use the representation of $A$ in Jordan normal form: 
\begin{align}  & \nn
A=UJU^{-1}, \quad J=\diag\{J_1, \, J_2, \ldots, J_N\}, \\
& \label{DD3}
 J_i= \la_i I_{n_i}+K_i, \quad K_i:=\begin{bmatrix} 0 & 1 & 0 & 0 &\ldots \\
0& 0& 1 & 0 & \ldots \\  \cdot & \cdot & \cdot & \cdot & \cdot
\end{bmatrix}.
\end{align}
The Jordan representation above yields the equality
\begin{align}& \label{DD4}
\E^{- \I t A}=U\diag\left\{\E^{- \I \la_1t }\sum_{k=0}^{n_1 -1}\frac{(-\I t K_i)^k}{k!}, \ldots, \E^{- \I \la_Nt }\sum_{k=0}^{n_N -1}\frac{(-\I t K_i)^k}{k!}\right\}
U^{-1}.
\end{align}
Taking into account formula \eqref{D5}, Corollary \ref{CyDD} and representation \eqref{DD4}, we see that 
the following asymptotics  is valid generically$:$
$$\|\psi(x,t)g\|=C_{\pm}(g)\E^{\tau_{\pm} t} |t|^{r_{\pm}} \big(1+O(1/t)\big) \quad {\mathrm{for}} \quad t\to \pm \infty,$$
where $\| \cdot \|$ is the norm in $L^2_h(0,\infty)$, $g\in \BC^h$, $\tau_+=\max_{1\leq i \leq N}\Im(\la_i)$, \\
$\tau_-=\min_{1\leq i \leq N}\Im(\la_i)$,  $\, \displaystyle{r_{\pm}=\max_{\Im(\la_i)=\tau_{\pm}}(n_i-1)}$.
\end{Rk}
We note that in a different way the Jordan structure of $A$ was used in \cite{ALS-KP} to study (and explain) an interesting multi-lump
phenomena discovered in \cite{ACTV}.
\paragraph{3.}
Using considerations similar to those in Paragraph 1 of this section, we construct GBDT for matrix Schr\"odinger equation \eqref{h1}.
Solution $w$ of system (\ref{bt1}) with the coefficients given by
(\ref{bt11}) can be written down in the block form: $w= \left[
\begin{array}{c} y \\  \hat{y} \end{array} \right]$ ($y, \hat y \, \in \BC^h$). Hence, we
rewrite (\ref{bt1}) as 
$$y^{\prime}(x, \lambda )= \hat{y}(x,
\lambda ),  \quad \hat{y}^{\prime}(x, \lambda )= - \lambda y(x,
\lambda )+u(x) y(x, \lambda ),
$$
that is, (\ref{h1}) is fulfilled. So
system (\ref{bt1}), (\ref{bt11}) is equivalent to the
Schr\"odinger equation (\ref{h1}). 
The following proposition is a
corollary of Theorem \ref{Tmbt}. 
\begin{Pn} \label{Pnbt}
Let a vector function $y(x, \lambda )$ satisfy the Schr\"odinger
equation \eqref{h1}, where the $h \times h$ potential $u=u^*$ is locally summable on $\BR$,  and 
assume that the conditions of Theorem \ref{Tmbt} hold.
Then, the vector function
\begin{equation} \label{bt13}
\widetilde{y}(x, \lambda )= [I_{h} \hspace{1em}0] \widetilde{w}(x, \lambda
), \quad \widetilde w(x, \lambda ):=w_A(x, \lambda )w(x, \lambda ), \end{equation}
with $w=\begin{bmatrix} y \\ y^{\prime} \end{bmatrix}$, satisfies the matrix
Schr\"odinger equation
\begin{equation} \label{bt14}
-\widetilde{y}^{\prime \prime }(x, \lambda )+ \widetilde{u}(x) \widetilde{
y}(x, \lambda )= \lambda \widetilde{ y}(x, \lambda ),
\end{equation}
where $\wt u =\wt u^*$ is given by the formula \eqref{bt15}.
\end{Pn}
\begin{proof}. According to Theorem \ref{Tmbt}, we have $\widetilde{w}^{\prime}= \widetilde{G} \widetilde{w}$. 
We rewrite this equation in terms of the blocks $\wt y$ and $\breve y:=[0 \quad I_{h} ] \widetilde{w}$ of $\widetilde{w}$:
$$ \wt y^{\prime}=-X_{22}\wt y +\breve y, \quad \breve y^{\prime}=-\la \wt y+ (u+X_{12}+X_{21})\wt y+X_{22} \breve y.$$
Differentiating $\wt y$ in the first equation above and using the second equation, we obtain $\widetilde{y}^{\prime \prime }
=-\la \wt y +(u+X_{12}+X_{21}+X_{22}^2-X_{22}^{\prime})\wt y.$
Now, using \eqref{h9}, we derive
\begin{align} & \label{h7}
\widetilde{y}^{\prime \prime }=-\la \wt y +\big(u+2\big(X_{12}+X_{21}+X_{22}^2\big)\big)\wt y.
\end{align}
Relation \eqref{bt14}  is immediate from \eqref{bt15} and \eqref{h7}.
\end{proof}
Instead of the Schr\"odinger equation (\ref{bt14}) one can talk
about Schr\"odinger operator ${\wt \cll}=-\frac{d}{dx^2} +{\wt u}$
with a properly defined domain.
%%%%%%%%%%%%%%%%%%%%%%%%%%%%%%%%%%%
\paragraph{4.} If $\Pi$ and $S$ are known explicitly, then representations \eqref{D5} and \eqref{bt13} provide explicit solutions of Schr\"odinger systems \eqref{I2'}
and \eqref{bt14}, respectively. In particular, $\Pi$ and $S$ are easily constructed in the case $u \equiv 0$ (see \cite{GKS3}). For this purpose
we partition $\Pi$ into $n \times h$ blocks: $\Pi=\begin{bmatrix}\Lam_1 & \Lam_2 \end{bmatrix}$. Then, the first equation in \eqref{h2} takes (for $u\equiv 0$) the form
\begin{align} & \label{D7} 
\Lam_1^{\prime}=A\Lam_2, \quad \Lam_2^{\prime}=-\Lam_1.
\end{align}
\begin{Rk} \label{RkExplPi}
When $u\equiv 0$, then $\Pi(x)$ in Theorems \ref{Tmbt} and  \ref{TmCS} and in Proposition \ref{Pnbt} is given by the formulas $\Pi(x)=\begin{bmatrix}\Lam_1(x) & \Lam_2(x) \end{bmatrix}$
and
\begin{align} & \label{D8} 
 \begin{bmatrix}\Lam_1(x) \\ \Lam_2(x) \end{bmatrix}=\E^{x\cla}\begin{bmatrix}\Lam_1(0) \\ \Lam_2(0) \end{bmatrix},
\quad \cla:=\begin{bmatrix}0 & A \\  -I_h & 0\end{bmatrix},
\end{align}
which is immediate from \eqref{D7}. According to \eqref{h2}, the matrix function $S(x)$ is given by the formula
\begin{align} & \label{D9} 
S(x)=S(0)+\int_0^x \Lam_2(\eta)\Lam_2(\eta)^*d\eta .
\end{align}
Recall that we know $($that is, we choose$)$ parameter matrices $A$, $S(0)$ and $\Pi(0)$ or, equivalently, $A$, $S(0)$ and $\Lam_k(0)$ $(k=1,2)$
which determine GBDT-transformation.
\end{Rk}
%%%%%%%%%%%%%%%%%%%%%%%%%%%%%%%%%%%%%%%%%%%
%%%%%%%%%%%%%%%%%%%%%%%%%%%%%%%%%%%%%%%%%%%%
\section{Discrete dynamical Schr\"odinger system} \label{DscrSchrEq}
\setcounter{equation}{0}
\paragraph{1.}
GBDT (generalized B\"acklund-Darboux transformation) was applied to important linear
and nonlinear discrete systems in \cite{FKKS, FKRS, KS, SaA07, SaSaR}.
In particular, discrete canonical systems and non-Abelian Toda lattices were studied
in \cite{SaA07}. Jacobi matrices corresponding to explicit solutions of matrix Toda lattices
were considered in \cite[Appendix]{SaA07}.  
Using some   modification of the results from \cite[Appendix]{SaA07}, we 
 construct here explicit solutions of discrete dynamical Schr\"odinger systems. We present also direct proofs 
 of the corresponding modified results from \cite[Appendix]{SaA07}, whereas in \cite[Appendix]{SaA07}
 several essential facts are proved indirectly (via the theory of discrete canonical systems developed in the previous sections
 of \cite{SaA07}) and some details of the proofs are  omitted.

We start  with introducing  generalized B\"acklund-Darboux transformation  (GBDT)
of block Jacobi matrices.
Suppose that the sets of $h \times h$ matrices $\{C(k)
\}_{k > 0}$ and $\{Q(k) \}_{k > 0}$ such that
\begin{equation} \label{A-3}
C(k)Q(k)^*=Q(k)C(k) \quad (k>0) , \quad C(k)>0 \quad (k > 0),
\end{equation}
are given.  The corresponding initial Jacobi matrix is introduced by the relations
\begin{align} & \label{A1-}
 \clj=\left[
\begin{array}{cccccc}
 b_1 &  a_1 & 0 & 0 & 0 & \ldots \\  c_2 &  b_2 & 
a_2 & 0 & 0 & \ldots \\  0 &  c_3 &  b_3 &  a_3 & 0 &
\ldots
\\ \ldots & \ldots & \ldots & \ldots & \ldots & \ldots
\end{array}
\right], \quad  c_k=  a_{k-1}^*,
\\ & \label{A2-}
 a_k=-\I  C(k)^{- \frac{1}{2}} C(k+1)^{\frac{1}{2}}, \quad
 b_k=  C(k)^{- \frac{1}{2}}  Q(k) 
C(k)^{\frac{1}{2}}, 
\end{align}
where $k>0$ and (according to \eqref{A-3} and \eqref{A2-}) $b_k=b_k^*$.

Recall that GBDT is determined by three parameter matrices. Thus, we fix $n>0$, two $n \times n$ parameter matrices $A$ and
$S_0>0$ and an $n \times m $ $(m=2h)$ parameter matrix $ \Pi_0$ such that
\begin{equation} \label{A-2}
AS_0-S_0A^*=\I \Pi_0  j \Pi_0^*, \quad j=\begin{bmatrix}0 & I_h \\ I_h & 0 \end{bmatrix}.
\end{equation}
Everywhere in this section $j$
is given by the second equality in \eqref{A-2}.
Introduce matrices $ \Pi_k$ and $S_k$ for $k>0$
by the recursions
\begin{equation} \label{A-1}
 \Pi_k=  \Pi_{k-1} \xi(k)^{-1}-\I A   \Pi_{k-1} { P },
\quad S_k=S_{k-1}+  \Pi_{k-1} \zeta(k)   \Pi_{k-1}^*,
\end{equation}
where
\begin{equation} \label{A}
\xi(k)= \left[
\begin{array}{lr}
- \I Q(k) & C(k) \\  C(k)^{-1} & 0
\end{array}
\right], \quad \zeta(k)=\left[
\begin{array}{lr}
0 & 0 \\ 0 & C(k)^{-1}
\end{array}
\right], \quad  P=\begin{bmatrix}0 & 0 \\ 0 &  I_h  \end{bmatrix}.
\end{equation}
 The following properties easily follow from \eqref{A-3} and \eqref{A}: $jPj=I_m-P$,
\begin{align} & \label{J6}
\xi(k)j\xi(k)^*=\xi(k)^*j\xi(k)=j,  \quad P\xi(k)j=\zeta(k), \quad PjP=0.
\end{align}
Therefore, taking adjoints of both parts of the first equality in \eqref{A-1} (and multiplying the result by $\I^k j$)
we obtain an equivalent to this equality relation
\begin{align} & \label{J14}
 \I^k j\Pi_k^*= \I \xi(k)\big(\I^{k-1}j \Pi_{k-1}^*\big) -  (I_m-P)   \big(\I^{k-1}j \Pi_{k-1}^*\big)A^* \qquad (m=2h). 
\end{align}
\begin{Rk}\label{MTL} Setting in \eqref{J14} $ \I^k j\Pi_k^*=W(k)$ and $A^*=z$,
we obtain an auxiliary linear system $($10.1.9$)$ from \cite{SaL30} for the matrix Toda chain,
which  explains the choice of  the equation on $\Pi_k$ in \eqref{A-1}. Namely, we see that this
equation is a generalized auxiliary system for Toda chain with the generalized
eigenvalue $A$.
\end{Rk}
Since $S_0>0$ and $C(k)>0$, the second equality in \eqref{A-1} yields $S_k>0$ for $k \geq 0$. Setting 
\begin{align} & \label{D10} 
X(k)= \{ X_{ip}(k)
\}_{i,p=1}^2= \Pi_k^* S_k^{-1}  \Pi_k \quad (k \geq 0),
\end{align}
we define the transformed matrices $\wt C(k)$ and $\wt Q(k)$ via relations
\begin{align} & \label{A0}
\wt C(k)=C(k)+X_{22}(k-1) \quad (k>0), 
\\ & \label{A0'}
\wt Q(k)=Q(k)+\I \big(X_{21}(k-1)-X_{21}(k)\big) \quad (k>0).
\end{align}
Clearly $X(k)\geq 0$ for $k\geq 0$, and so  $\wt C(k)>0$ for $k > 0$.
Then, the transformation (GBDT) $\wt \clj$ of the block Jacobi matrix $ \clj$ is defined by the equalities
\begin{align} & \label{A1}
\wt \clj=\left[
\begin{array}{cccccc}
\wt b_1 & \wt a_1 & 0 & 0 & 0 & \ldots \\ \wt c_2 & \wt b_2 & \wt
a_2 & 0 & 0 & \ldots \\  0 & \wt c_3 & \wt b_3 & \wt a_3 & 0 &
\ldots
\\ \ldots & \ldots & \ldots & \ldots & \ldots & \ldots
\end{array}
\right], \quad \wt c_k= \wt a_{k-1}^*,
\\ & \label{A2}
\wt a_k=-\I \wt C(k)^{- \frac{1}{2}}\wt C(k+1)^{\frac{1}{2}}, \quad
\wt b_k= \wt C(k)^{- \frac{1}{2}} \wt Q(k) \wt
C(k)^{\frac{1}{2}}.
\end{align}
We note that formulas (\ref{A1}) and (\ref{A2}) coincide (after removal of tildes) with the  formulas (\ref{A1-}) and (\ref{A2-})
which define $ \clj$.

According to \cite[Appendix]{SaA07},
we have $\wt b_k=\wt b_k^*$.  
Slightly modifying the proof of 
\cite[Theorem A.1]{SaA07}, one may derive that under condition
\begin{align} & \label{J0}
\begin{bmatrix}I_h &0\end{bmatrix}\Pi_0^*S_0^{-1}=0
\end{align}
we have
\begin{align} & \label{J1}
\wt \clj Y=YA \qquad \big(Y=\{y_k\}_{k \geq 1}, \quad y_{k}:= \begin{bmatrix} 0\ \wt C(k)^{-1/2}\end{bmatrix}\Pi_{k-1}^*S_{k-1}^{-1}\big).
\end{align}

\begin{Tm} \label{TmA0}
Suppose that Jacobi matrix $\wt \clj$ is given by the formulas
\eqref{A1} and \eqref{A2}, that relations \eqref{A-3}, \eqref{A-2} and \eqref{J0}
are valid, and that the matrices $\wt C(k)$ and $\wt Q(k)$ in \eqref{A2}
are given by \eqref{A-1}--\eqref{A0'}.

Then the block vector function
\begin{align} & \label{J2}
\Psi(t)=Y\E^{-\I t A},
\end{align}
where $Y$ is introduced in \eqref{J1}, satisfies the
discrete dynamical Schr\"odinger system $\I \left(\frac{\p}{\p t}\Psi\right)(t)=\wt \clj \Psi(t)$.
\end{Tm}
Theorem \ref{TmA0} is immediate from \eqref{J1} and it remains to prove \eqref{J1}.
More precisely, we prove the following theorem.
\begin{Tm} \label{TmA} Suppose that the relations \eqref{A-3} and \eqref{A-2} $($where $S_0>0)$
are valid. Then the matrices $\wt C(k)$ and $\wt Q(k)$ given by \eqref{A-1}--\eqref{A0'}
are well-defined and satisfy relations
\begin{equation} \label{J3}
\wt C(k) \wt Q(k)^*=\wt Q(k) \wt C(k) \quad (k>0) , \quad \wt C(k)>0 \quad (k > 0).
\end{equation}
We also have $\wt b_k=\wt b_k^*$ for the matrices $\wt b_k$ in \eqref{A2}.  

If, in addition, \eqref{J0} holds, then the matrix $\wt \clj$ of the form \eqref{A1}, \eqref{A2} satisfies \eqref{J1}.
\end{Tm}
\begin{proof}. Step 1. 
Taking into account the inequalities $S_0>0$, $C(k)>0$ and relations \eqref{A-1}, \eqref{A0}, we explained already 
that $S_k>0$ $(k\geq 0)$ and that  $\wt C(k)>0$ $(k >0)$. Therefore, the matrices $\wt Q(k)$ and $\wt C(k)$ are well-defined,
the inequality for $\wt C(k)$ in \eqref{J3} is valid, $\wt C(k)$ is invertible, and $\wt \clj$ is also well-defined. 

Using \eqref{A-2} we show by induction that the matrix identity
\begin{equation} \label{A-2+}
AS_k-S_kA^*=\I \Pi_k  j \Pi_k^*
\end{equation}
holds for all $k\geq 0$.
Namely, let us assume that the identity 
$$AS_{k-1}-S_{k-1}A^*=\I \Pi_{k-1}  j \Pi_{k-1}^*$$ 
is valid for some $k>0$. Then, in view of 
the second equality in \eqref{A-1} 
we have
\begin{align} \nn
AS_k-S_kA^*&=AS_{k-1}-S_{k-1}A^*+A\Pi_{k-1}\zeta(k)\Pi_{k-1}^*-\Pi_{k-1}\zeta(k)\Pi_{k-1}^*A^*
\\ & \label{J4}
=\I \Pi_{k-1}j\Pi_{k-1}^*+A\Pi_{k-1}\zeta(k)\Pi_{k-1}^*-\Pi_{k-1}\zeta(k)\Pi_{k-1}^*A^*.
\end{align}
On the other hand, the first equality in \eqref{A-1} and relations \eqref{J6} imply that
\begin{align} & \label{J5}
\I \Pi_{k}j\Pi_{k}^*==\I \Pi_{k-1}j\Pi_{k-1}^*+A\Pi_{k-1}\zeta(k)\Pi_{k-1}^*-\Pi_{k-1}\zeta(k)\Pi_{k-1}^*A^*.
\end{align}
(Here we used also the equality $\xi(k)^{-1}=j\xi(k)^*j$, which is immediate from \eqref{J6}.)
Comparing \eqref{J4} and \eqref{J5} we obtain \eqref{A-2+}.

Step 2. Next, we prove the equality
\begin{align} & \label{J7}
\Pi_k^*S_k^{-1}=\I P \Pi_{k-1}^*S_{k-1}^{-1}A+ j \wt \xi(k) j \Pi_{k-1}^*S_{k-1}^{-1}, 
\\ & \label{J8}
\wt \xi(k):=\begin{bmatrix} -\I \wt Q(k)  & \wt C(k)\\ \breve C(k) & 0
\end{bmatrix}, \quad \breve C(k)=C(k)^{-1}-X_{11}(k).
\end{align}
Indeed, taking into account the second relation in \eqref{A-1}, 
the equality $$j\xi(k)^* P=\zeta(k),$$ which is immediate from 
 \eqref{J6}, and the equalities  
 $$ j(I_m-P)j=P, \quad P(I_m-P)=0,$$ 
 we derive
\begin{align} \nn
S_k^{-1}-S_{k-1}^{-1} &=-S_k^{-1}\Pi_{k-1}\zeta(k)\Pi_{k-1}^*S_{k-1}^{-1}
\\ & \label{J9}
=
-S_k^{-1}\big(\Pi_{k-1}j\xi(k)^*j -\I A \Pi_{k-1} P\big) (I_m-P)j\Pi_{k-1}^*S_{k-1}^{-1}.
\end{align}
In view of the first relation in \eqref{A-1} and the equality $j\xi(k)^*j=\xi(k)^{-1}$, we rewrite \eqref{J9}
in the form
\begin{align} & \label{J10}
S_k^{-1}-S_{k-1}^{-1} =-S_k^{-1}\Pi_{k}(I_m-P)j\Pi_{k-1}^*S_{k-1}^{-1}.
\end{align}
Multiplying both sides of  \eqref{J10} by $\Pi_k^*$ from the left and using again the first relation in \eqref{A-1},
we see that
\begin{align} \label{J11}
\Pi_k^*S_k^{-1}&=\Pi_k^*S_{k-1}^{-1} -X(k)(I_m-P)j\Pi_{k-1}^*S_{k-1}^{-1}
\\ & \nn
=
j\xi(k)j \Pi_{k-1}^*S_{k-1}^{-1}+\I P\Pi_{k-1}^*A^*S_{k-1}^{-1}-X(k)(I_m-P)j\Pi_{k-1}^*S_{k-1}^{-1}.
\end{align}
Substituting (into \eqref{A-2+}) $k-1$ instead of $k$, we rewrite
the result in the form 
\begin{align}& \label{J12}
A^*S_{k-1}^{-1}= S_{k-1}^{-1}A- \I S_{k-1}^{-1} \Pi_{k-1}j\Pi_{k-1}^*S_{k-1}^{-1}.
\end{align}
After substitution of \eqref{J12} into \eqref{J11}, we obtain
\begin{align}  \nn
\Pi_k^*S_k^{-1}=&\I P\Pi_{k-1}^*S_{k-1}^{-1}A+j\xi(k)j \Pi_{k-1}^*S_{k-1}^{-1}+PX(k-1)j\Pi_{k-1}^*S_{k-1}^{-1}
\\ & \label{J13}
-X(k)(I_m-P)j\Pi_{k-1}^*S_{k-1}^{-1}.
\end{align}
Equality \eqref{J7} follows from \eqref{J13}.

Step 3. Recall that $\xi(k)$ is $j$-unitary, that is, the relation $\xi(k)j\xi(k)^*=j$
(or, equivalently, $\xi(k)^*j\xi(k)=j$) in \eqref{J6} holds. It is important to show that the transformed matrix $\wt \xi(k)$ introduced
in \eqref{J8} is $j$-unitary as well.  Taking into account \eqref{A0} and \eqref{A0'} we rewrite  \eqref{J8} in the form
\begin{align}& \label{J16}
\wt \xi(k)=\xi(k)-jX(k)(I_m-P)+jPX(k-1).
\end{align}

Assuming that $\det A\not=0$, we prove the equality
\begin{align}& \label{J17}
\wt \xi(k)=\breve w(k)\xi(k)\breve w(k-1)^{-1}, \quad \breve w(k):=I_m-\I j \Pi_k^*S_k^{-1}A^{-1}\Pi_k.
\end{align}
We note that $\breve w(k)=w_A(k,0)$, where $w_A(k, \la)=I_m-\I j \Pi_k^*S_k^{-1}(A-\la I_n)^{-1}\Pi_k$
(with matrices $A, \Pi_k, S_k$ satisfying \eqref{A-2+}) is the transfer matrix function in Lev Sakhnovich form.
(Compare with $w_A$ in \eqref{bt5}.) According to \cite{SaL1} (see also \cite[p. 24]{SaSaR}) we have
$w_A(k,\la)jw_A(k, \ov \la)^*=j$, and so the matrices $\breve w(k)$ ($k \geq 0$) are $j$-unitary.
Thus, \eqref{J17} implies that $\wt \xi(k)$ is $j$-unitary. If $\det A=0$, we have $\det(A-\la I_n)\not=0$
for smal values $\la=\ov \la$ and approximate $\wt \xi(k)$ with the $j$-unitary matrices $\xi_{\la}(k)$
corresponding to the GBDT-generating triples $A-\la I_n, \, S_0, \, \Pi_0$. Hence, the equality \eqref{J17} for 
the case $\det A\not= 0$ yields the $j$-unitarity property
\begin{align}& \label{J18}
\wt \xi(k)j\wt \xi(k)^*=\wt \xi(k)^*j \wt \xi(k)=j
\end{align}
without restriction on $\det A$. It remains to show that \eqref{J17} is valid.

Indeed, let us rewrite \eqref{J17} in the form
\begin{align}& \label{J19}
\wt \xi(k)(I_m-\I j \Pi_{k-1}^*S_{k-1}^{-1}A^{-1}\Pi_{k-1})=(I_m-\I j \Pi_k^*S_k^{-1}A^{-1}\Pi_k) \xi(k).
\end{align}
Using the first relation in \eqref{A-1}, we rewrite \eqref{J19} in another equivalent form
\begin{align}\nn
\wt \xi(k)-\I\wt \xi(k) j \Pi_{k-1}^*S_{k-1}^{-1}A^{-1}\Pi_{k-1}=&\xi(k)-\I j \Pi_k^*S_k^{-1}A^{-1}\Pi_{k-1}
\\ & \label{J20}
-  j \Pi_k^*S_k^{-1}\Pi_{k-1}P\xi(k).
\end{align}
Substituting the expression for $\Pi_k^*S_k^{-1}$ from \eqref{J7} into the second right hand term
in \eqref{J20}, we obtain
\begin{align}\nn
\wt \xi(k)&-\I\wt \xi(k) j \Pi_{k-1}^*S_{k-1}^{-1}A^{-1}\Pi_{k-1}=\xi(k)+ jP \Pi_{k-1}^*S_{k-1}^{-1}\Pi_{k-1}
\\ & \label{J21}
-\I \wt \xi(k) j\Pi_{k-}^*S_{k-1}^{-1}A^{-1}\Pi_{k-1}-  j \Pi_k^*S_k^{-1}\Pi_{k-1}P\xi(k).
\end{align}
Therefore, using \eqref{J16} and canceling similar terms we derive
\begin{align}\nn &
-jX(k)(I_m-P)+jPX(k-1)=jP \Pi_{k-1}^*S_{k-1}^{-1}\Pi_{k-1}-  j \Pi_k^*S_k^{-1}\Pi_{k-1}P\xi(k).
\end{align}
Recalling that by definition $X(k-1)= \Pi_{k-1}^*S_{k-1}^{-1}\Pi_{k-1}$ we further simplify the equality
above, and so the following relation:
\begin{align}\label{J22} &
-j \Pi_k^*S_k^{-1}\Pi_{k}(I_m-P)=-  j \Pi_k^*S_k^{-1}\Pi_{k-1}P\xi(k)
\end{align}
is equivalent to \eqref{J17}. In view of \eqref{A-1} and \eqref{A}, we have
$$\Pi_{k}(I_m-P)=\Pi_{k-1}P\xi(k),$$ which proves \eqref{J22}. Thus, \eqref{J17}
is proved as well, and hence \eqref{J18} holds.

In particular, \eqref{J18} yields
\begin{align}\label{J23} &
(I_m-P)\wt \xi(k)j\wt \xi(k)^*(I_m-P)=0, \quad P\wt \xi(k)j\wt \xi(k)^*(I_m-P)=I_h.
\end{align}
According to \eqref{J8}, formula \eqref{J23} is equivalent to the equalities
\begin{align}\label{J24} &
\wt C(k)\wt Q(k)^*-\wt Q(k)\wt C(k)=0, \quad \breve C(k)=C(k)^{-1}.
\end{align}
Therefore, the first equality in \eqref{J3} is valid, and we also rewrite \eqref{J8}  in the form
\begin{align}
 & \label{J8'}
\wt \xi(k)=\begin{bmatrix} -\I \wt Q(k)  & \wt C(k)\\ \wt C(k)^{-1} & 0
\end{bmatrix}.
\end{align}
Comparing \eqref{J8'} and the first equality in \eqref{A}, we see that the representations of $\wt \xi(k)$
and $\xi(k)$ differ only by tildes in the notations.

The equality  $\wt b_k= \wt b_k^*$ (for $\wt b_k$ given by  (\ref{A2})) is immediate
from the first relation in \eqref{J3}. 

%%%%%%%%%%%%%%%%%%%%%%%%%%%%%%%%
Step 4. Finally, let us prove \eqref{J1}.
Using equalities  (\ref{A2}) and the definition of $Y$ in (\ref{J1}),
we obtain
\begin{align}
 & \label{J25}
\wt a_k y_{k+1}=-\I \wt C(k)^{-\frac{1}{2}}
\begin{bmatrix} 0   & I_h
\end{bmatrix}\Pi_{k}^* S_{k}^{-1}
\\  & \label{J26}
\wt b_k y_{k}=\wt C(k)^{-\frac{1}{2}}\wt Q(k)
\begin{bmatrix} 0   & I_h
\end{bmatrix}\Pi_{k-1}^* S_{k-1}^{-1}.
\end{align}
Relations \eqref{J7}, \eqref{J8} and \eqref{J25}, \eqref{J26} imply that
\begin{align}
 \nn
\wt a_k y_{k+1}+\wt b_k y_{k}= &\wt C(k)^{-\frac{1}{2}}
\begin{bmatrix} 0   & I_h
\end{bmatrix}\Pi_{k-1}^* S_{k-1}^{-1}A
\\ & \label{J27}
-\I \wt C(k)^{\frac{1}{2}}
\begin{bmatrix} I_h   & 0
\end{bmatrix}\Pi_{k-1}^* S_{k-1}^{-1}.
\end{align}
In particular, taking into account \eqref{J0}, we derive
\begin{align}
 \label{J28}
\wt b_1 y_{1}+ \wt a_1 y_{2} =\wt C(1)^{-\frac{1}{2}}
\begin{bmatrix} 0   & I_h
\end{bmatrix}\Pi_{0}^* S_{0}^{-1}A=y_1 A.
\end{align}

From \eqref{J7} and \eqref{J8'} we see that
\begin{align}
 \label{J29}
\begin{bmatrix} I_h   & 0
\end{bmatrix}\Pi_{k-1}^* S_{k-1}^{-1}= \wt C(k-1)^{-1}
\begin{bmatrix} 0   & I_h
\end{bmatrix}\Pi_{k-2}^* S_{k-2}^{-1} \quad (k>1).
\end{align}
According to \eqref{A2} and \eqref{J29} we have
\begin{align}
\nn
a_{k-1}^*y_{k-1}&=\I \wt C(k)^{\frac{1}{2}}\wt C(k-1)^{-1}
\begin{bmatrix} 0   & I_h
\end{bmatrix}\Pi_{k-2}^* S_{k-2}^{-1}
\\  \label{J30} &
=
\I \wt C(k)^{\frac{1}{2}}\begin{bmatrix} I_h   & 0
\end{bmatrix}\Pi_{k-1}^* S_{k-1}^{-1}.
\end{align}
Now, formulas \eqref{J27} and \eqref{J30}
yield for $k>1$ that
\begin{align}
& \label{J31}
a_{k-1}^*y_{k-1}+\wt b_k y_{k}+\wt a_k y_{k+1}= &\wt C(k)^{-\frac{1}{2}}
\begin{bmatrix} 0   & I_h
\end{bmatrix}\Pi_{k-1}^* S_{k-1}^{-1}A=y_kA.
\end{align}
Equalities  \eqref{J28} and \eqref{J31} imply \eqref{J1}.
\end{proof}

\noindent{\bf Acknowledgments.}
 {The research of A.L. Sakhnovich   was supported by the
Austrian Science Fund (FWF) under Grant  No. P29177.}

\newpage

\begin{flushright}
B. Fritzsche,\\
Fakult\"at f\"ur Mathematik und Informatik,  Universit\"at Leipzig, \\
 Augustusplatz 10,  D-04009 Leipzig, Germany, \\
 e-mail: {\tt Bernd.Fritzsche@math.uni-leipzig.de}
 
\vspace{0.5em}

B. Kirstein, \\ 
Fakult\"at f\"ur Mathematik und Informatik,  Universit\"at Leipzig, \\
 Augustusplatz 10,  D-04009 Leipzig, Germany, \\
 e-mail: {\tt Bernd.Kirstein@math.uni-leipzig.de}
 
\vspace{0.5em}

I.Ya. Roitberg, \\
 e-mail: {\tt  	innaroitberg@gmail.com}

\vspace{2em} 

A.L. Sakhnovich,\\
Faculty of Mathematics,
University
of
Vienna, \\
Oskar-Morgenstern-Platz 1, A-1090 Vienna,
Austria, \\
e-mail: {\tt oleksandr.sakhnovych@univie.ac.at}

\end{flushright}


\begin{thebibliography}{AGKS}
\bibitem{ACTV}
M.J.~Ablowitz, S.~Chakravarty, A.D.~Trubatch  and J.~Villarroel, A novel class of solutions of the non-stationary
Schr\"odinger and the Kadomtsev--Petviashvili I equations, {\it Phys. Lett. A} {\bf 267} (2000), 132--146.


\bibitem{AMRa}
 S.~Avdonin, V.~Mikhaylov  and K.~Ramdani,  Reconstructing the potential for the one-dimensional Schr\"odinger equation from boundary measurements, 
\textit{IMA J. Math. Control Inform.} \textbf{31} (2014), 137--150.

 \bibitem{BeMi}
 M.~Belishev and V.~Mikhailov,  Inverse problem for a one-dimensional dynamical Dirac system (BC-method), \textit{Inverse Problems} 
 \textbf{30} (2014), 125013.
 
\bibitem{BoPP} 
M. Boiti, F. Pempinelli and A.K. Pogrebkov, 
On the extended resolvent of the nonstationary Schr\"odinger operator for a Darboux transformed potential,
\textit{J. Phys. A}  {\bf 39}:8 (2006),  1877--1898. 

\bibitem{Ci}
J.L. Cieslinski, 
{Algebraic construction of the Darboux matrix revisited}, 
 \textit{J. Phys. A} \textbf{42}(40) (2009), 404003, 40 pp.

  \bibitem{D}
P.A. Deift,  
{Applications of a commutation formula}. \textit{ Duke Math. J.} \textbf{45} (1978), 267--310.

\bibitem{EKMT}
 I.E.~Egorova,   E.A.~Kopylova,  V.A.~Marchenko and G.~Teschl, 
Dispersion estimates for one-dimensional Schr\"odinger and Klein--Gordon equations revisited, \textit{Russian Math. Surveys} \textbf{ 71}:3 (2016), 391--415.

\bibitem{EKT}
  I.E.~Egorova,   E.A.~Kopylova and G.~Teschl,  Dispersion estimates for one-dimensional discrete Schr\"odinger and wave equations, \textit{J.~Spectr. Theory}
  \textbf{5}:4  (2015), 663--696.

\bibitem{FKKS} 
B. Fritzsche, M.A. Kaashoek, B. Kirstein and  A.L. Sakhnovich,  Skewselfadjoint Dirac systems with rational rectangular Weyl functions: explicit solutions of direct and inverse problems and integrable wave equations,  \textit{Math. Nachr.} \textbf{289}:14-15 (2016), 1792--1819.  

\bibitem{FKRS}
B. Fritzsche, B. Kirstein, I. Roitberg and A.L. Sakhnovich, 
Weyl matrix functions and  inverse problems for discrete Dirac type self-adjoint system: explicit and general 
solutions. \textit{Oper. Matrices} \textbf{2} (2008), 201--231.

\bibitem{Ge}
F. Gesztesy,  
{A complete spectral characterization of the double commutation method}. \textit{J. Funct. Anal.}  \textbf{117}:2 (1993), 401--446.

\bibitem{GeT}
F. Gesztesy and G. Teschl,  
On the double commutation method, \textit{Proc. Amer. Math. Soc.} \textbf{124} (1996), 1831--1840.

\bibitem{GKS3}
 I. Gohberg, M.A.  Kaashoek and A.L.  Sakhnovich, Sturm--Liouville systems with rational Weyl functions: explicit formulas and applications,  \textit{Integral Equations Operator Theory}  \textbf{30}:3 (1998), 338--377. 

\bibitem{Gu}
C.H.~Gu, H.~Hu and  Z.~Zhou, 
\textit{Darboux transformations in integrable systems. Theory and their applications to geometry}, Mathematical Physics Studies \textbf{26},  Springer, Dordrecht, 2005.

\bibitem{JK}
A. Jensen and T. Kato, Spectral properties of Schr\"odinger operators and time-decay of
the wave functions, {\it Duke Math. J.} {\bf 46} (1979), 583--611.

\bibitem{KS}
M.A.  Kaashoek and A.L.  Sakhnovich,  Discrete skew self-adjoint canonical system and the isotropic Heisenberg magnet model,  \textit{J. Funct. Anal.} \textbf{228}:1 (2005),  207--233.

\bibitem{KT} 
E.A. Kopylova and G. Teschl,  Dispersion estimates for one-dimensional discrete Dirac equations,  \textit{J. Math. Anal. Appl.} \textbf{434}:1 (2016), 191--208. 

\bibitem{KoSaTe}
A. Kostenko,  A. Sakhnovich  and  G. Teschl, 
Commutation methods for Schr\"odinger operators with strongly singular potentials, \textit{Math. Nachr.} \textbf{285} (2012), 392--410.
 
\bibitem{Mar}
V.A. Marchenko,
\textit{ Nonlinear equations and operator algebras}, Mathematics and Its Applications (Soviet Series) \textbf{17},  D. Reidel, Dordrecht etc., 1988.


\bibitem{MS}
V.B. Matveev  and M.A.  Salle,    \textit{Darboux transformations and solitons}, Springer, Berlin,  1991.
 

\bibitem{PrdO}  
 R.A. Prado and  C.R.  de Oliveira,  Sparse 1D discrete Dirac operators I: Quantum transport,  \textit{J. Math. Anal. Appl.} \textbf{385}:2  (2012),  947--960.
 
\bibitem{SaA2}
A.L. Sakhnovich,    Dressing procedure for solutions of
nonlinear equations and the method of operator identities, \textit{Inverse problems}
\textbf{10} (1994), 699--710.


\bibitem{SaA3}
A.L. Sakhnovich,    Generalized B\"acklund--Darboux transformation: spectral properties and nonlinear equations, \textit{J. Math. Anal. Appl.} \textbf{262} (2001), 274--306.


 \bibitem{ALS-PT}
A.L. Sakhnovich,     Non-Hermitian matrix Schr\"odinger equation: B\"acklund-Darboux transformation, Weyl functions and PT-symmetry,  
\textit{J. Phys. A} \textbf{36}:28  (2003),  7789--7802.

\bibitem{ALS-KP}
A.L. Sakhnovich, {  Matrix Kadomtsev-Petviashvili equation:
matrix identities and explicit non-singular solutions},  {\it J. Phys.
A} {\bf 36} (2003), 5023--5033.


\bibitem{SaA07}
A.L. Sakhnovich,    
Discrete canonical system and non-Abelian Toda lattice: B\"acklund--Darboux transformation and Weyl functions,  \textit{  Math. Nachr.} \textbf{280}:5-6 (2007), 1--23.

\bibitem{SaA15}
 A.L. Sakhnovich,     Dynamical and spectral Dirac systems: response function and inverse problems,  \textit{J. Math. Phys.} \textbf{56}:11  (2015),  Paper 112702, 13 pp.
 
\bibitem{SaA17}   
A.L. Sakhnovich, {Dynamics of electrons and explicit
solutions of Dirac--Weyl systems}, \textit{J. Phys. A: Math. Theor.} {\bf 50} (2017), Paper 115201.
 

\bibitem{SaA16}
A.L. Sakhnovich,    Dynamical canonical systems and their explicit solutions,
 \textit{Discrete and Continuous Dynamical Systems Series A} {\bf 37}:3 (2017), 551--561.


%\bibitem{SaAHams}
%A.~L. Sakhnovich,    
%Hamiltonian systems and Sturm-Liouville equations: Darboux transformation and applications, \textit{Integral Equations Operator Theory}
%DOI 10.1007/s00020-017-2385-7.
 
 
\bibitem{SaSaR}
A.L. Sakhnovich,   L.A.  Sakhnovich    and I. Roitberg,  \textit{Inverse Problems and Nonlinear Evolution Equations. 
 Solutions, Darboux Matrices and Weyl--Titchmarsh Functions},  {De Gruyter Studies in Mathematics} \textbf{47},   De Gruyter, Berlin, 2013.
 
 \bibitem{SaL1}
L.A. Sakhnovich,     On  the  factorization  of  the  transfer matrix function, \textit{ Sov. Math. Dokl.} \textbf{17} (1976), 203--207.

\bibitem{SaL30}
L.A.  Sakhnovich, 
\textit{ Interpolation theory and its applications},  Mathematics and its Applications \textbf{428}, Kluwer,  Dordrecht,  1997.



\bibitem{SaL2}
L.A. Sakhnovich,   \textit{Spectral theory of canonical differential systems, method of operator identities}, Operator Theory Adv. Appl. \textbf{107}, Birkh\"auser, Basel, 1999.


\bibitem{T0}
G. Teschl, 
\textit{Jacobi operators and completely integrable nonlinear lattices},  Mathematical Surveys and Monographs \textbf{72},  Amer. Math. Soc., Providence, 2000.

\bibitem{T1}
 G. Teschl, \textit{Ordinary differential equations and dynamical systems},  Graduate Studies in Mathematics \textbf{140},  Amer.  Math. Soc., Providence, RI, 2012.

\end{thebibliography}
\end{document}